\documentclass[10pt]{amsart}
\usepackage{amsfonts,amstext,amssymb}
\renewcommand\newsymbol[5]{%
\DeclareMathSymbol#1{#3}%
   {\ifcase #2\or AMSa\or AMSb\fi}{"#4#5}}
\newsymbol\bbox 1004
\newsymbol\beperk 1216
\def\concat{{}^\frown}
\renewenvironment{proof}{\par\medskip\noindent\textbf{Proof}}
                       {\hfill$\bbox$\par\medskip}
\newtheorem{definition}{Definition}
\newtheorem{theorem}[definition]{Theorem}
\newtheorem{lemma}[definition]{Lemma}
\newtheorem{corollary}[definition]{Corollary}
\newtheorem{remark}[definition]{Remark}
\newcommand\stem{\operatorname{stem}}
\newcommand\dom{\operatorname{dom}}
\newcommand\cf{\operatorname{cf}}
\newcommand{\powerinfront}[1]{{}^#1 2}
\renewcommand\int{\operatorname{int}}
\begin{document}
\begin{abstract}
In this paper we show that, when we iteratively add Sacks reals to
a model of $ZFC$ we have for every two reals in the extension a
continuous function defined in the ground model that maps one of
the reals onto the other.
\end{abstract}
\title{A small transitive family of continuous functions on the Cantor set}
\author{K.P. Hart and B.J. van der Steeg}
\maketitle

%
%
\section{Introduction}
In \cite{D} Dow gave a proof of the Rudin-Shelah theorem about the
existence of $2^\mathfrak{c}$ points in $\beta\mathbb{N}$ that are
Rudin-Keisler incomparable. The proof actually shows that whenever
a family $\mathcal{F}$ of $\mathfrak{c}$ continuous self-maps of
$\beta\mathbb{N}$ (or $\mathbb{N}^*$) are given there is a set $S$
of $2^\mathfrak{c}$ many $\mathcal{F}$-independent points in
$\beta\mathbb{N}$ (or $\mathbb{N}^*$). This suggests that we
measure the complexity of a space $X$ by the cardinal number
$\mathfrak{t}\mathfrak{f} (X)$ defined as the minimum cardinality
of a set $\mathcal{F}$ of continuous self maps such that for all
$x,y\in X$ there is $f\in\mathcal{F}$ such that $f(x)=y$ or
$f(y)=x$. Let us call such an $\mathcal{F}$ transitive. Thus Dow's
proof shows $\mathfrak{t}\mathfrak{f}(\beta\mathbb{N})$,
$\mathfrak{t}\mathfrak{f}(\mathbb{N}^*)\geq\mathfrak{c}^+$.

We investigate $\mathfrak{t}\mathfrak{f}(C)$, where $C$ denotes
the Cantor set. Van Mill observed that
$\mathfrak{t}\mathfrak{f}(C)\geq\aleph_1$; a slight extension of
his argument shows that $\operatorname{MA(countable)}$ implies
$\mathfrak{t}\mathfrak{f}(C)=\mathfrak{c}$. Our main result states
that in the Sacks model the continuous functions on the Cantor set
that are coded in the ground model form a transitive set. Thus we
get the consistency of
$\mathfrak{t}\mathfrak{f}(C)=\aleph_1<\aleph_2=\mathfrak{c}$.

The gap between $\mathfrak{t}\mathfrak{f}(C)$ and $\mathfrak{c}$
cannot be arbitrarily wide, because Hajnal's free set lemma
implies that for any space $X$ one has $\vert X\vert\leq
\mathfrak{t}\mathfrak{f}(X)^+$.

In \cite{M} Miller showed that it is consistent with $ZFC$ that
for every set of reals of size continuum there is a continuous map
from that set onto the closed unit interval. In fact he showed
that the iterated perfect set model of Baumgartner and Laver (see
\cite{BL}) is such a model, and noted that the continuous map can
even be coded in the ground model.

Here we will show that in the iterated perfect set model, for
every two reals $x$ and~$y$ there exists a continuous function
with code in the ground model that maps $x$ onto~$y$ or $y$
onto~$x$.
\begin{definition}
By a transitive set of functions $\mathcal{F}$ we mean a set of
continuous functions such that for every two reals $x$ and~$y$
there exists an element $f\in\mathcal{F}$ such that $f(x)=y$ or
$f(y)=x$ holds.

Let us also define the cardinal number $\mathfrak{t}\mathfrak{f}$
by
$$ \mathfrak{t}\mathfrak{f}=\min\{\vert\mathcal{F}\vert
:\mathcal{F}\ \text{is a transitive set of functions}\},\
\text{i.e.}\
\mathfrak{t}\mathfrak{f}=\mathfrak{t}\mathfrak{f}(C).$$
\end{definition}
The paper is organized as follows: in section~2 we prove some
simple facts on $\mathfrak{t}\mathfrak{f}$, the minimal size of
transitive sets of functions. We also state and prove the main
theorem of this paper in section~2, using theorems proved later on
in section~3. As a corollary to the main theorem we have the
consistency of $\mathfrak{t}\mathfrak{f}<\mathfrak{c}$ with $ZFC$.
Finally in section~4 we will make a remark on the effect on
$\mathfrak{t}\mathfrak{f}$ when we add $\kappa$ many Sacks reals
side-by-side to a model of $ZFC+CH$.
%
%
\section{Notation and Preliminaries}
For the rest of this paper let $V$ be a model of $ZFC$. We will
use the same notations and definitions as Baumgartner and Laver in
\cite{BL}, so for any ordinal $\alpha$ we let $\mathbb{P}_\alpha$
denote the poset that iteratively adds $\alpha$ Sacks reals to the
model $V$, using countable support. Let $\mathbb{P}_1=\mathbb{P}$,
where $\mathbb{P}$ denotes the 'normal' Sacks poset for the
addition of one Sacks real.

Let $G_\alpha$ be $\mathbb{P}_{\alpha}$-generic over~$V$, we
define $V_\alpha$ by $V_\alpha=V[G_\alpha]$ for every ordinal
$\alpha$. Note that if $\beta<\alpha$ we have that
$G_\alpha\beperk\beta$ is a $\mathbb{P}_\beta$-generic subset
over~$V$. If we denote the $(\alpha+1)$-th added Sacks real by
$s_\alpha$ then we can also write $V_\alpha=V[\langle
s_\beta:\beta<\alpha\rangle]$.

Assuming $V\models CH$, the proof of the following facts can be
found in \cite{BL}:
\begin{enumerate}
\label{fact}
\item
    Forcing with $\mathbb{P}_\alpha$ does not
    collapse cardinals.
\item
    $V_{\omega_2}$ is a model of $ZFC+2^{\aleph_0}=\aleph_2$.
\item
    Let $\dot{\mathbb{P}}_\beta$ denote the result of
    defining $\mathbb{P}_\beta$ in $V_\alpha$. Then for any
    $\alpha,\beta\geq 1$, $\Vdash_\alpha
    ``\mathbb{P}_{\alpha,\alpha+\beta}\ \text{is isomorphic to}\
    \dot{\mathbb{P}}_\beta "$.
\end{enumerate}
We will now prove some facts on the cardinal
$\mathfrak{t}\mathfrak{f}$. The first is Van Mill's observation
alluded to above.
\begin{theorem}
\label{fgealeph1} $\mathfrak{t}\mathfrak{f}\ge\aleph_1$.
\end{theorem}
\begin{proof}
Suppose $\mathcal{F}$ is a countable set of functions. Let $A_f$
denote the set $\{x:\int(f^{-1}(x))\not=\emptyset\}$ for every
$f\in\mathcal{F}$. Every $A_f$ is at most countable because
$2^\omega$ is separable. So choose an $x$ in
$2^\omega\setminus\bigcup_{f\in\mathcal{F}}A_f$, then we know that
for every $f\in\mathcal{F}$ the set $f^{-1}(x)$ is nowhere dense
in $2^\omega$. For such an $x$ the set $\{f^{-1}(x):
f\in\mathcal{F}\}$ is countable. Because the set $\{f(x):
f\in\mathcal{F}\}$ is also countable the Baire category theorem
tells us that  the set $2^\omega\setminus
\bigcup_{f\in\mathcal{F}}(\{f(x)\}\cup f^{-1}(x))$ is nonempty,
thus showing that $\mathcal{F}$ is not transitive.
\end{proof}
\begin{theorem}
\label{f<=c<=f+}
$\mathfrak{t}\mathfrak{f}\leq\mathfrak{c}\leq\mathfrak{t}\mathfrak{f}^+$.
\end{theorem}
\begin{remark}
\label{remarkcov} The proof of theorem~\ref{fgealeph1} shows that
$\mathfrak{t}\mathfrak{f}$ is at least the minimum number of
nowhere dense sets needed to cover~$C$. Then
theorem~\ref{f<=c<=f+} and $\operatorname{MA(countable)}$ imply
$\mathfrak{t}\mathfrak{f}=\mathfrak{c}$.
\end{remark}
The second inequality is a consequence of the following lemma. The
proof of this lemma can be found in \cite{W}.

For this we need some more notation. Let $S$ be an arbitrary set.
By a \emph{set mapping on $S$} we mean a function $f$ mapping $S$
into the power set of~$S$. The set map is said to be of
\emph{order $\lambda$} if $\lambda$ is the least cardinal such
that $\vert f(x)\vert<\lambda$ for each $x$ in~$S$. A subset $S'$
of~$S$ is said to be \emph{free for $f$} if for every $x\in S'$ we
have $f(x)\cap S'\subset\{x\}$.
\begin{lemma}[Free set lemma]
Let $S$ be a set with $\vert S\vert=\kappa$ and $f$ a set map on
$S$ of order $\lambda$ where $\lambda <\kappa$. Then there is a
free set of size $\kappa$ for~$f$.
\end{lemma}
\begin{proof}[Theorem~\ref{f<=c<=f+}]
The proof of the first inequality is easy. We simply have to
observe that the set of all constant functions on the reals is a
transitive set of functions.

Now for the second inequality. Striving for a contradiction
suppose that $\mathfrak{c}\geq\mathfrak{t}\mathfrak{f}^{++}$. Let
$\mathcal{F}$ be a transitive set of functions such that
$\vert\mathcal{F}\vert=\mathfrak{t}\mathfrak{f}$. We define a set
map $F$ on the reals by $ F(x)=\{f(x): f\in\mathcal{F}\}$ for
every $x\in 2^\omega$. Because $\vert
F(x)\vert\leq\mathfrak{t}\mathfrak{f}$, this set map $F$ is of
order $\mathfrak{t}\mathfrak{f}^+$, which is less than
$\mathfrak{c}$. According to the free set lemma there exists a set
$X\subset 2^\omega$ such that $\vert X\vert=\mathfrak{c}$ and for
every $x\in X$ we have $F(x)\cap X\subset\{x\}$. This is a
contradiction, because every two reals in $X$ provide a counter
example of $\mathcal{F}$ being a transitive set.
\end{proof}
Closed subsets of the Cantor set can be coded by sub trees of
$\powerinfront{{<\omega}}$, as follows: if $A$ is closed then let
$T_A=\{x\beperk n:\ x\in A,\ n\in\omega\}$; one can recover $A$
from $T_A$ by observing that $A=\{x\in\powerinfront{\omega}:\
\forall n\in\omega,\ x\beperk n\in T_A\}$.

When we say that a closed set $A$ is \emph{coded in the ground
model} we mean that $T_A$ belongs to the ground model.

We shall always construct a continuous function $f$ between closed
sets $A$ and~$B$ by specifying an order-preserving map $\phi$ from
$T_A'$ to $T_B$, where $T_A'$ denotes the set of splitting nodes
of $T_A$. Once $\phi$ is found one defines $f$ by
$$f(x)=``\text{the path through}\ T_B\ \text{determined by the
restriction of}\ \phi\ \text{to}\ \{x\beperk n:\ n\in\omega\}".$$
We say that $f$ is coded in the ground model if $\phi$ belongs to
$V$. In what follows we shall denote the map $\phi$ by~$f$ as
well.

Let us define the set $\mathcal{G}$ (in any $V_\alpha$) by
$$ \mathcal{G}=\{f:f\ \text{is a continuous function with code
in}\ V\}$$
Now we can explicitly state the main theorem of this paper.
Section~\ref{proofMainth} is completely devoted to the proof of
this theorem by parts, so we will prove the theorem here and refer
to the needed theorems proved in that section.
\begin{theorem}[Main Theorem]
\label{thmain} The set $\mathcal{G}$ is transitive in $V_\alpha$
for every ordinal $\alpha$.
\end{theorem}
\begin{proof}
We will show by transfinite induction that $\mathcal{G}$ is a
transitive set in $V_\alpha$ for~all~$\alpha$. For $\alpha=0$ this
is obvious. Suppose the theorem is true for all $\beta<\alpha$.
Let $x$ and~$y$ be reals in $V_\alpha$.

If $\alpha$ is  a successor ordinal, $\alpha=\beta+1$, then we use
theorem~\ref{thopv} in the case that at least one of $x$ and~$y$
is not in $V_\beta$ to show that there exist a continuous function
$f$ defined in~$V$ (so $f\in\mathcal{G}$) such that in $V_\alpha$
we have $f(x)=y$ or $f(y)=x$.

Since we are forcing with countable support and because reals are
countable objects, there are no new reals added by
$\mathbb{P}_\alpha$ for $\cf(\alpha)>\aleph_0$. So if $\alpha$ is
a limit ordinal we only have to consider the case where
$\cf(\alpha)=\aleph_0$ and at least one of $x,y$ is not in
$\bigcup_{\beta<\alpha}V_\beta$. Then we use theorem~\ref{thlim}
to show the existence of an continuous function~$f$ defined in~$V$
such that in $V_\alpha$ $f(x)=y$ or $f(y)=x$ holds.
\end{proof}
As is well-known, if $V\models CH$ then
$V_{\omega_2}\models\mathfrak{c}=\aleph_2$. This enables us to
show that $\mathfrak{t}\mathfrak{f}<\mathfrak{c}$ is consistent.
\begin{corollary}
If $V\models CH$ then
$V_{\omega_2}\models\mathfrak{t}\mathfrak{f}<\mathfrak{c}$.
\end{corollary}

In this paper we shall repeatedly use the fact that any
homeomorphism~$h$ between two closed nowhere dense subsets of the
Cantor set can be extended to a homeomorphism of the Cantor set
onto itself (see \cite{KR}). Furthermore it is straight forward to
extend a continuous function between to closed nowhere dense
(disjoint) subsets of the Cantor set to a continuous self map of
the Cantor set.

Because we can make sure that the subsets of the Cantor set that
define the added reals $x$ and~$y$ are nowhere dense and closed,
when we show that there exists a homeomorphism (or a continuous
function) $f$ mapping of one of these sets onto the other, in such
a way that in the extension $x$ is mapped onto~$y$ or vice versa,
we actually have shown that there exists a self map of the Cantor
that is a homeomorphism (continuous function) mapping, in the
extension, $x$ onto~$y$ or $y$ onto~$x$.
%
%
\section{The continuous functions with code in the ground model~$V$
        form a transitive set in $V_\alpha$}
\label{proofMainth} In this section we prove that for every
$\alpha$ and any new real $x$ in the Baumgartner and Laver model
$V_{\alpha}$ (i.e. $x\in
V_\alpha\setminus\bigcup_{\beta<\alpha}V_\beta$) and $y$ any real
in $V_\alpha$ there exists a function $f$ defined in the ground
model $V$ such that in $V_\alpha$ the equation $f(x)=y$ holds.

We make the following definition. For any
$\sigma\in\powerinfront{{<\omega}}$ we let $l(\sigma)\in\omega$
denote the \emph{length of $\sigma$}. So for every
$\sigma\in\powerinfront{{<\omega}}$ we have
$\sigma\in\powerinfront{{l(\sigma)}}$.
%
%
To show how we construct our continuous maps we reprove the
familiar fact that Sacks reals are minimal, see~\cite{J}.
\begin{lemma}
\label{lem01} Suppose $x$ is a real in $V[G]\setminus V$, where
$G$ is a $\mathbb{P}$-generic filter over~$V$, and that
$p\in\mathbb{P}$ is such that $p\Vdash ``\dot{x}\not\in V"$. Then
there exists a $q\geq p$ and a homeomorphism~$f$ defined in~$V$
such that $q\Vdash ``f(\dot{s})=\dot{x}"$. Here $\dot{s}$ denotes
the name of the added Sacks real.
\end{lemma}
\begin{proof}
We will construct a fusion sequence $\{\langle p_i,n_i\rangle:
i\in\omega\}$ such that each $p_{i+1}$ will know all the first $i$
splitting nodes of every branch of the perfect tree $p_i$ and
$(p_{i+1},n_{i+1})>(p_i,n_i)$ for every $i$.

Because $p$ forces that $\dot{x}$ is a new real, there exists an
element $u_\emptyset\in\powerinfront{{<\omega}}$ with maximal
length $m_\emptyset$, such that $p\Vdash ``\dot{x}\beperk
m_\emptyset=u_\emptyset"$ and $p$ does not decide
$\dot{x}(m_\emptyset)$. There exist $p_{\langle
0\rangle},p_{\langle 1\rangle}\geq p_0$ such that $p_{\langle
k\rangle}\Vdash ``\dot{x}(m_\emptyset)=k"$ for $k\in\{0,1\}$.
Without loss of generality the stems of $p_{\langle 0\rangle}$ and
$p_{\langle 1\rangle}$ are incompatible. Let
$n_0=\min\{n\in\omega: p_{\langle 0\rangle}\beperk n\not=
p_{\langle 1\rangle}\beperk n\}$ and let $p_0$ denote the element
$p_{\langle 0\rangle}\cup p_{\langle 1\rangle}$.

Now assume we have
$p_i=\bigcup\{p_\sigma:\sigma\in\powerinfront{{i+1}}\}$. Consider
$\tau\in\powerinfront{{i+1}}$, we have an element
$u_\tau\in\powerinfront{{<\omega}}$ of maximal length $m_\tau$
such that $p_\tau\Vdash ``\dot{x}\beperk m_\tau=u_\tau"$. There
exist $p_{\tau\concat{0}}$, $p_{\tau\concat{1}}\geq p_\tau$ such
that $p_{\tau\concat{k}}\Vdash ``\dot{x}(m_\tau)=k"$ for
$k\in\{0,1\}$. Again without loss of generality the stems of
$p_{\tau\concat{0}}$ and $p_{\tau\concat{1}}$ are incompatible.
Let $n_\tau$ denote the integer $\min\{n\in\omega:
p_{\tau\concat{0}}\beperk n\not=p_{\tau\concat{0}}\beperk n\}$ and
$n_{i+1}=\max\{n_\sigma:\sigma\in\powerinfront{{i+1}}\}$. We let
$p_{i+1}$ denote the element
$\bigcup\{p_\sigma:\sigma\in\powerinfront{{i+2}}\}$. Now the
induction step is completed, because $p_{i+1}$ knows all the first
$i+1$ splitting nodes of every branch in $p_i$ and
$(p_{i+1},n_{i+1})>(p_i,n_i)$ for every $i\in\omega$.

We define the function~$f$ by
$$ f^{-1}([u_\sigma])\supset[\stem(p_\sigma)]\ \text{for}\
\sigma\in\powerinfront{{<\omega}}.$$
As $\stem(p_\sigma)$ is a finite approximation of the added Sacks
real $\dot{s}$, we have by the construction of our $p_\sigma$ for
$\sigma\in\powerinfront{{<\omega}}$ and the function $f$ that
$p_\sigma\Vdash ``f(\dot{s})\in [u_\sigma]"$ for every
$\sigma\in\powerinfront{{<\omega}}$. And so the fusion $q$ of the
sequence $\{\langle p_i,n_i\rangle: i\in\omega\}$ forces that in
the extension $V[G]$ the equality $f(s)=x$ holds. This $f$, being
a continuous bijection between two Cantor sets, is (of course) a
homeomorphism.
\end{proof}
\begin{remark}
\label{remarktree}In the lemma we have also defined a map $\phi$
from the finite sub trees of the fusion $q$ onto the finite sub
trees of $T=\bigcup_{\sigma\in\powerinfront{{<\omega}}}u_\sigma$
which induces our homeomorphism. We have $\phi(q)=T$ and
$$\phi([q\beperk\sigma])=\bigcup\{u_\tau:\sigma\subset\tau\
\text{and}\ \tau\in\powerinfront{{<\omega}} \}.$$
We note that $[T]$ is the set of all the possible interpretations
of $\dot{x}$ in~$V[G]$ and that $T$ depends on $\phi$ and~$q$
only. In theorem~\ref{thopv} we will use this interpretation of
the previous lemma.
\end{remark}
As a warming up exercise we prove the following.
\begin{theorem}
\label{thstap0} The set $\mathcal{G}$ is transitive in~$V_1$.
\end{theorem}
\begin{proof}
Suppose $x$ and~$y$ are two reals of $V_1(=V[s_0])$. We consider
two cases.\\
{\bf Case 1:} $x$ is a real in $V$. The constant function
$c_x=\{\langle y,x\rangle: y\ \text{a real in}\ V_1\}$ is a
continuous function defined in $V$, thus a member of
$\mathcal{G}$, and in~$V_1$ it maps $y$ onto~$x$.
\par\medskip\noindent {\bf Case 2:} both $x$ and~$y$ are reals not
in~$V$. Let $p\in\mathbb{P}$ be a witness of this, so $p\Vdash
``\dot{x},\dot{y}\not\in V"$. According to lemma~\ref{lem01} there
exists a $q\geq p$ and a homeomorphism~$f$ defined in~$V$ such
that $q\Vdash ``f(\dot{s}_0)=\dot{x}"$, where $\dot{s}_0$ denotes
the added Sacks real. If we apply the lemma again we get an $r\geq
q$ and a homeomorphism~$g$ defined in~$V$ such that $q\Vdash
``g(\dot{s}_0)=\dot{y}"$. But now we have that $r\Vdash ``(g\circ
f^{-1})(\dot{x})=\dot{y}"$ and we see that $g\circ f^{-1}$ is the
element of $\mathcal{G}$ we are looking for.
\end{proof}
%
%
\begin{theorem} \label{thopv} For $\alpha$ an ordinal and
$x$ and~$y$ reals in~$V_{\alpha+1}$ such that $x\not\in V_\alpha$
there exists an $f\in\mathcal{G}$ such that in~$V_{\alpha+1}$
$f(x)=y$ holds.

Moreover if also $y\not\in V_\alpha$ then $f$ can be chosen to be
a homeomorphism.
\end{theorem}
\begin{proof}
This is an immediate consequence of the lemmas~\ref{lem1} and
\ref{lem2}.
\end{proof}
We make the following definitions. For $p\in\mathbb{P}$ and
$s\in\powerinfront{{<\omega}}$ we let $p_s$ denote the sub-tree
$\{t\in p:s\subseteq t\ \text{or}\ t\subseteq s\}$ of~$p$. Of
course $p_s$ is a perfect tree if and only if
$s\in\powerinfront{{<\omega}}\cap p$. To generalize this to
$\mathbb{P}_\alpha$, suppose $p$ is an element of
$\mathbb{P}_\alpha$, $F$ is a finite subset of $\dom(p)$ and
$n\in\omega$, we say that a function $\tau:F\to\powerinfront{n}$
\emph{is consistent with} $p$ if the following holds for every
$\beta\in F$:
$$(p\beperk\tau)\beperk\beta\Vdash_\beta``\tau(\beta)\in
p(\beta)".$$ So we have for every $\beta\in F$ that
$(p\beperk\tau)\beperk\beta\Vdash_\beta``(p(\beta))_{\tau(\beta)}\
\text{is a perfect tree}"$.

Furthermore let us suppose that $F$ and~$H$ are two sets such that
$F\subset H$, and $n$ and~$m$ are two integers such that $m<n$, if
$\tau$ is a function mapping $F$ into~$\powerinfront{m}$ then we
say that a function $\sigma: H\rightarrow\powerinfront{n}$
\emph{extends} the function $\tau$ if for every $i\in F$ we have
$\sigma(i)\beperk m=\tau(i)$.

For later use we will prove the following:
\begin{lemma}
\label{smelt} Let $p\in\mathbb{P}_\alpha$,
$F\in[\dom(p)]^{<\omega}$ and $n\in\omega$. Suppose
$\tau:F\to\powerinfront{n}$ is consistent with~$p$ then for every
$r\geq p\beperk\tau$ there exists a $q\geq p$ such that
$q\beperk\tau=r$ and $q\beperk\beta\Vdash_\beta
``(p(\beta))_s=(q(\beta))_s\ \text{for every}\
s\in\powerinfront{n}$ such that $s\not= \tau(\beta)"$ for every
$\beta\in F$.
\end{lemma}
\begin{proof}
Define the element $q\in\mathbb{P}_\alpha$ as follows for
$\beta<\alpha$:
\begin{displaymath}
q\beperk\beta\Vdash_\beta ``q(\beta)=\left\{
    \begin{array}{ll}
    r(\beta) & \beta\not\in F\\
    r(\beta)\cup\{(p(\beta))_s:s\in\powerinfront{n}\cap p(\beta)\
    \text{such that}\ s\not=\tau(\beta) &  \beta\in F".
\end{array}
\right.
\end{displaymath}
In this way we strengthen the tree $p(\beta)$ above $\tau(\beta)$
keeping the rest of the perfect tree intact (according to $F$
anyway).
\end{proof}
We need the following lemma to make sure that the maps we will
construct in the lemmas~\ref{lem1} and \ref{lem2} are well-defined
and continuous.
\begin{lemma}
\label{lemhlpopv} Let $p\in\mathbb{P}_{\alpha+1}$. Suppose
$F,H\in[\dom(p)]^{<\omega}$ are such that $F\subset H$ and
$m,n\in\omega$ are such that $m<n$. If $\tau:
F\rightarrow\powerinfront{m}$ is consistent with $p$, $N$ is an
integer and $T$ is a finite tree such that
$$(p\beperk\tau)\beperk\alpha\Vdash ``p(\alpha)\cap
\powerinfront{{\leq N}}=T",$$
then there exist a $(q,j)>_H(p\beperk\tau,n)$ and an $M>N$ such
that for every $\sigma:H\rightarrow\powerinfront{n}$ extending
$\tau$, if $\sigma$ is consistent with~$q$, then there exists
$T_\sigma$ such that $q\beperk\sigma\Vdash
``q(\alpha)\cap\powerinfront{{\leq M}}=T_\sigma"$. Also
$\vert(T_\sigma)_t\cap\powerinfront{M}\vert\geq 2$ for every $t\in
T$ and $[T_\sigma]\cap[T_\varsigma]=\emptyset$ whenever $\sigma$
and~$\varsigma$ are distinct and consistent with~$q$.
\end{lemma}
\begin{proof}
Let $\Sigma_\tau$ denote the set of all
$\sigma:H\rightarrow\powerinfront{n}$ extending~$\tau$. Because
$p(\alpha)$ is a perfect tree there exists a
$\mathbb{P}_\alpha$-name $\dot{M}$ such that for every $t\in T$ we
have
$$(p\beperk\tau)\beperk\alpha\Vdash ``\vert (p(\alpha))_t\cap
\powerinfront{{\dot{M}}}\vert\geq 2\vert\Sigma_\tau\vert".$$
According to lemma 2.3 of \cite{BL} there exists a
$(q^\dagger,j^\dagger)>_H((p\beperk\tau)\beperk\alpha,n)$ such
that if $\sigma\in\Sigma_\tau$ is consistent with~$q^\dagger$ we
have an $M_\sigma$ such that $q^\dagger\beperk\sigma\Vdash
``\dot{M}=M_\sigma"$. Put $M=\max\{M_\sigma:\sigma\in\Sigma_\tau\
\text{consistent with}\ q^\dagger\}$. We have $q^\dagger\Vdash
``\vert (p(\alpha))_t\cap\powerinfront{M}\vert\geq
2\vert\Sigma_\tau\vert"$ for every $t\in T$.

Enumerate $\{\sigma\in\Sigma_\tau:\sigma\ \text{consistent with}\
q^\dagger\}$ as $\{\sigma_k:k<K\}$. Let $r\geq
q^\dagger\beperk\sigma_0$ be such that $r\Vdash
``p(\alpha)\cap\powerinfront{{\leq M}}=S_{\sigma_0}"$, where
$S_{\sigma_0}$ is such that $\vert (S_{\sigma_0})_t\cap
\powerinfront{M}\vert\geq 2\vert\Sigma_\tau\vert$ for every $t\in
T$. Use lemma~\ref{smelt} to find a $q_0\geq q^\dagger$ such that
$q_0\beperk\sigma_0=r$.

We continue this procedure with all the $\sigma_k\in\Sigma_\tau$.
So if $\sigma_k$ is consistent with~$q_{k-1}$ we find an $r\geq
q_{k-1}\beperk\sigma_k$ such that $r\Vdash
``p(\alpha)\cap\powerinfront{{\leq M}}=S_{\sigma_k}"$, and also
that $\vert (S_{\sigma_k})_t\cap\powerinfront{M}\vert\geq
2\vert\Sigma_\tau\vert$ for every $t\in T$. And we use
lemma~\ref{smelt} to define $q_k\geq q_{k-1}$ such that
$q_k\beperk\sigma_k=r$. If $\sigma_k$ is not consistent with
$q_{k-1}$ we choose $q_k=q_{k-1}$.

We now have for every $\sigma\in\Sigma_\tau$ consistent with
$q_{K-1}$ a finite tree $S_\sigma\subset \powerinfront{{\leq M}}$
extending the tree $T$ such that every branch in $T$ has (at
least) $2\vert\Sigma_\tau\vert$ different extensions in
$S_\sigma\cap\powerinfront{M}$ and $q_{K-1}\beperk\sigma\Vdash
``p(\alpha)\cap\powerinfront{{\leq M}}=S_\sigma"$.

As $q_{K-1}$ forces that, for each $y\in T$ the size of the set
$p(\alpha))_t\cap\powerinfront{M}$ is at least
$2\vert\Sigma_\tau\vert$ we can find for $\sigma\in\Sigma_\tau$
consistent with $q_{K-1}$ a sub tree $T_\sigma$ of $S_\sigma$ such
that $\vert(T_\sigma)_t\cap\powerinfront{M}\vert\geq 2$ and
whenever $\sigma$ and $\varsigma$ are distinct and consistent with
$q_{K-1}$ we have $[T_\sigma]\cap[T_\varsigma]=\emptyset$.

Define $q\in\mathbb{P}_{\alpha+1}$ such that $q\beperk\alpha=
q_{K-1}$ and choose $q(\alpha)$ such that for every consistent
$\sigma\in\Sigma_\tau$ we have $q\beperk\sigma\Vdash
``q(\alpha)=p(\alpha)\cap[T_\sigma]"$. If we let $j$ be equal to
$\max\{j^\dagger,M\}$ the proof is complete.
\end{proof}
\begin{lemma}
\label{lem1} Given an ordinal~$\alpha$, a
$p\in\mathbb{P}_{\alpha+1}$ and $\mathbb{P}_{\alpha+1}$-names
$\dot{x}$ and $\dot{y}$ such that $p\Vdash``\dot{x}\not\in
V_\alpha\ \text{and}\ \dot{y}\in V_\alpha"$ then there exists a
continuous function~$f$ defined in~$V$ and a $q\geq p$ such that
$q\Vdash``f(\dot{x})=\dot{y}"$.
\end{lemma}
\begin{proof}
By remark~\ref{remarktree} we know that there is an $r\geq
p\beperk\alpha$ and there exist $\mathbb{P}_{\alpha+1}$ names
$\dot{\phi}$ for a map on the finite sub trees of $p(\alpha)$ and
$\dot{T}$ for a perfect tree such that $r\Vdash
``\dot{\phi}(p(\alpha))=\dot{T}"$. Without loss of generality we
assume that $p\beperk\alpha=r$.

Let us construct a  fusion sequence $\{\langle
p_i,n_i,F_i\rangle:i\in\omega\}$. Let $p_0=p_1=p$, $n_0=n_1=0$,
$F_0=\emptyset$ and choose $F_1\in[\dom(p)]^{<\omega}$ in such a
way that we are building a fusion sequence.

Suppose we have constructed the sequence up to~$i$, let us
construct the next element of the fusion sequence. We let
$\{\tau_k:k<K\}$ denote all
$\tau:F_{i-1}\rightarrow\powerinfront{{n_{i-1}}}$ consistent
with~$p_i$. If we choose in lemma~\ref{lemhlpopv} $\tau=\tau_0$,
$F=F_{i-1}$ and $m=n_{i-1}$ we get a
$(q_0,m_0)>_{F_i}(p_i\beperk\tau_0,n_i)$ such that for every
$\sigma:F_i\to\powerinfront{{\leq n_i}}$ extending~$\tau_0$,
consistent with~$q_0$, we have a finite sub tree
$T_\sigma\subset\powerinfront{{\leq M(\tau_0)}}$
($M(\tau_0)\in\omega$ follows from lemma~\ref{lemhlpopv}) of
$p_i(\alpha)=p(\alpha)$ such that
\begin{enumerate}
\item
    $T_\sigma$ is an extension of $T_{\tau_0}$,
\item
    \label{itemref}
    for every branch $t$ in~$T_{\tau_0}$ there exist at least two
    different branches of length $M(\tau_0)$ in $T_\sigma$
    extending~$t$,
\item
    if $\sigma$ and~$\varsigma$ are two distinct members of
    $\Sigma_{\tau_0}$ consistent with~$q_0$ we have $[T_\sigma]\cap
    [T_\varsigma]=\emptyset$.
\end{enumerate}
We choose $r_0\in\mathbb{P}_{\alpha+1}$ with lemma~\ref{smelt}
such that $r_0\geq q_0$ and $r_0\beperk\tau_0=q_0$.

We iteratively consider all the
$\tau:F_{i-1}\to\powerinfront{{n_{i-1}}}$. In the general case if
$\tau_k$ is consistent with $r_{k-1}$ then lemma~\ref{lemhlpopv}
gives us a $q_k$ and an $m_k\in\omega$ such that
$(q_k,m_k)>_{F_i}(r_{k-1}\beperk\tau_k, n_i)$. We choose $r_k$ in
the same way as above, using lemma~\ref{smelt} such that $r_k\geq
q_k$ and $r_k\beperk\tau_k=q_k$. If $\tau_k$ is inconsistent
with~$r_{k-1}$ then we choose $r_k=r_{k-1}$ and $m_k=m_{k-1}$.
After considering all the $\tau_k$'s we define $p_{i+1}=r_{K-1}$
and $n_{i+1}=\max\{m_k: k<K\}$. This ends the construction of the
next element of the fusion sequence.

For every $i<\omega$ if $\sigma:F_i\to\powerinfront{{n_i}}$ is
consistent with~$p_{i+1}$ and extends $\tau:
F_{i-1}\to\powerinfront{{n_{i-1}}}$ then
$$ p_{i+1}\beperk\sigma\Vdash ``p(\alpha)\cap\powerinfront{{\leq
M(\tau)}}=T_\sigma".$$
Considering our function $\dot{\phi}$, let us denote the finite
tree $\dot{\phi}(T_\sigma)$ by $S_\sigma$. We have
$$p_{i+1}\beperk\sigma\Vdash ``\dot{\phi}(T_\sigma)=S_\sigma".$$
When we are building the fusion sequence we can of course make
sure that the fusion determines~$\dot{y}$ as well. Suppose we have
that $p_i\beperk\tau_k\Vdash ``t_{\tau_k}\subset\dot{y}"$,
$t_{\tau_k}$ of length~$i+1$. With lemma~\ref{lemhlpopv} we can
choose $q_k$ strong enough such that for every
$\sigma\in\Sigma_{\tau_k}$ consistent with~$q_k$ we have a
$t_\sigma$ of length~$i+2$ such that $q_k\beperk\sigma\Vdash ``
t_\sigma\subset\dot{y}"$. So assume we have made sure this is the
case and let us define the function~$f$ in~$V$ by $f(b)=t_\sigma$
for every maximal branch $b\in S_\sigma$ for every
$\sigma:F_i\to\powerinfront{{n_i}}$ consistent with~$p_i$ for some
$i\in\omega$. The function~$f$ is well-defined by
lemma~\ref{lemhlpopv} and we have for every $i\in\omega$ and
$\sigma:F_i\to\powerinfront{{n_i}}$ consistent with~$p_i$ that
$p_i\beperk\sigma\Vdash ``f([S_\sigma])\subset [t_\sigma]"$ and
thus $q\Vdash ``f(\dot{x})=\dot{y}"$.
\end{proof}
\begin{lemma}
\label{lem2} Given an ordinal~$\alpha$, a
$p\in\mathbb{P}_{\alpha+1}$ and $\mathbb{P}_{\alpha+1}$-names
$\dot{x}$ and~$\dot{y}$ such that $p\Vdash
``\dot{x},\dot{y}\not\in V_\alpha"$ then there exists a
homeomorphism~$f$, with code in~$V$, and a $q\geq p$ such that
$q\Vdash `` f(\dot{x})=\dot{y}"$.
\end{lemma}
\begin{proof}
By applying remark~\ref{remarktree} twice we have an $r\geq p$ in
$\mathbb{P}_{\alpha+1}$ and $\mathbb{P}_{\alpha+1}$ names
$\dot{\phi}_x$, $\dot{\phi}_y$ and~$\dot{T}_x$, $\dot{T}_y$ for
maps and perfect trees respectively such that
$r\beperk\alpha\Vdash ``\dot{\phi}_x(p(\alpha))=\dot{T}_x\
\text{and}\ \dot{\phi}_y(p(\alpha))=\dot{T}_y"$. Without loss of
generality we can assume that $p\beperk\alpha=r$.

During the construction of possible finite sub trees
$(T_x)_\sigma$ for~$\dot{x}$, when constructing the fusion
sequence in the proof of lemma~\ref{lem1} we could of course at
the same time also have constructed a similar sequence of finite
sub trees $(T_y)_\sigma$ for~$\dot{y}$.

Without loss of generality we could also have made sure that in
the proof of lemma~\ref{lem1} item~\ref{itemref} is replaced by
\begin{enumerate}
\item[$2^\dag$]
    for every maximal branch~$t$ in~$T_{\tau_0}$ there are exactly two
    different branches of length $M(\tau_0)$ in $T_\sigma$
    extending~$t$,
\end{enumerate}
Following the proof of lemma~\ref{lem1} we have for every
$\sigma:F_i\to\powerinfront{{n_i}}$ consistent with~$p_{i+1}$
finite sub trees $S_\sigma^x$ and~$S_\sigma^y$ such that $$
p_{i+1}\beperk\sigma\Vdash
``\dot{\phi}_x((T_x)_\sigma)=S_\sigma^x\ \text{and}\
\dot{\phi}_y((T_y)_\sigma)=S_\sigma^y.$$

We are ready to define the homeomorphism~$f$ in~$V$ that maps $x$
onto~$y$ in the extension. Suppose
$\tau:F_i\to\powerinfront{{n_i}}$ and
$\sigma:F_{i+1}\to\powerinfront{{n_{i+1}}}$ such that $\sigma$
extends~$\tau$. Every maximal branch in $(T_x)_\tau$ corresponds
to exactly one maximal branch in~$(T_y)_\tau$. Let $f$ map the
splitting point in $(T_x)_\sigma$ above any maximal branch in
$(T_x)_\tau$ onto the splitting point in $(T_y)_\sigma$ above the
corresponding maximal branch in~$(T_y)_\tau$. The function~$f$
thus defined will be a continuous and one-to-one mapping between
two Cantor sets, so a homeomorphism. Furthermore the fusion~$q$
forces that $f$ maps $x$ onto~$y$ in the extension.
\end{proof}
%
%
%
\begin{lemma}
\label{lemhlplim} Suppose that $\alpha$ is a limit ordinal of
cofinality $\aleph_0$. Let $x$ be a real in $V_\alpha$ such that
$x\not\in\bigcup_{\beta<\alpha}V_\beta$, and let
$p\in\mathbb{P}_\alpha$ be a witness of this. Also let $F, H\in
[\dom(p)]^{<\omega}$ such that $F\subset H$ and let $n$ and~$m$ be
two integers such that $m<n$. If $\tau:
F\rightarrow\powerinfront{m}$ is consistent with~$p$, and
$u_\tau\in\powerinfront{{<\omega}}$ is such that
$$p\beperk\tau\Vdash ``u_\tau\subset\dot{x}", $$
then there exists a $(q,j)>_H(p\beperk\tau,n)$ such that for every
$\sigma :H\rightarrow \powerinfront{n}$ consistent with~$q$, we
have a $u_\sigma\in\powerinfront{{<\omega}}$ such that
$q\beperk\sigma\Vdash ``u_\sigma\subset\dot{x}"$; in addition we
have $l(u_\sigma)=l(u_\varsigma)$ and $u_\sigma\not=u_\varsigma$
whenever $\sigma$ and~$\varsigma$ are distinct and consistent
with~$q$.
\end{lemma}
Before we prove the lemma we need some more notation. We let
$\Vdash^*$ denote forcing in $V_\delta$ over
$\mathbb{P}_{\delta\alpha}$. Here we use again the same notation
as in \cite{BL} where for $\delta<\alpha$
$P_{\delta\alpha}=\{p\in\mathbb{P}_\alpha: \dom(p)\subset\{\xi:
\delta\leq\xi<\alpha\}\}$, and if $p\in\mathbb{P}_\alpha$ then
$p^\delta=p\setminus(p\beperk\delta)\in\mathbb{P}_{\delta\alpha}$.
The mapping which carries $p$ into $(p\beperk\delta,p^\delta)$ is
an isomorphism of $\mathbb{P}_\alpha$ onto a dense subset of
$\mathbb{P}_\delta\times\mathbb{P}_{\delta\alpha}$
(see~\cite{BL}).
\begin{proof}[Lemma~\ref{lemhlplim}]
Choose a $\delta$ such that $\max(H)<\delta<\alpha$. Let $\tau:
F\rightarrow \powerinfront{m}$ be consistent with~$p$ and let
$\Sigma_\tau$ denote all the $\tau$ extending functions $\sigma:
H\rightarrow\powerinfront{n}$.

Because $p$ forces that $x\not\in V_\delta$, there is an antichain
below $p^\delta$ of size $\vert\Sigma_\tau\vert$ such that all
these elements force different interpretations of $\dot{x}$ in the
extension. In other words there exist a sequence
$\{\dot{f}_\sigma:\sigma\in\Sigma_\tau\}$ of $\mathbb{P}_\delta$
names for elements of $\mathbb{P}_{\delta\alpha}$ and a sequence
$\{\dot{u}_\sigma:\sigma\in\Sigma_\tau\}$ of $\mathbb{P}_\delta$
names for elements of $\powerinfront{{<\omega}}$ such that for all
$\sigma\in\Sigma_\tau$ we have
\begin{equation}
\label{equ1} (p\beperk\tau)\beperk\delta\Vdash
``\dot{f}_\sigma\geq p^\delta\ \text{and}\ \dot{f}_\sigma\Vdash^*
``\dot{u}_\sigma\subset\dot{x}"",
\end{equation}
and if $\sigma$ and~$\varsigma$ are distinct then
\begin{equation}
\label{equ2}(p\beperk\tau)\beperk\delta\Vdash
``l(\dot{u}_\sigma)=l(\dot{u}_\varsigma)\ \text{and}\
\dot{u}_\sigma\not=\dot{u}_\varsigma".
\end{equation}
Repeatedly using lemma~2.3 of \cite{BL} we see that there exist a
$(q^\dagger,j)>_H ((p\beperk\tau)\beperk\delta,n)$ and sequences
$\{f_\sigma:\sigma\in\Sigma_\tau\}$,
$\{u_\sigma:\sigma\in\Sigma_\tau\}\subset\powerinfront{i}$ for
some integer~$i$ such that for every $\sigma\in\Sigma_\tau$ we
have
\begin{equation}
\label{equ3}q^\dagger\Vdash_\delta ``\dot{f}_\sigma=f_\sigma\
\text{and}\ \dot{u}_\sigma=u_\sigma".
\end{equation}

Now let $q$ denote the element of~$\mathbb{P}_\alpha$ such that
$q\beperk\delta =q^\dagger$, and
$(q\beperk\sigma)\beperk\delta\Vdash ``q^\delta=f_\sigma"$ for
every $\sigma\in\Sigma_\tau$ consistent with~$q^\dagger$. This
completes the proof.
\end{proof}
\begin{theorem}
\label{thlim} For $\alpha$ a limit ordinal of cofinality
$\aleph_0$ and $x$ and~$y$ reals in~$V_\alpha$ such that
$x\not\in\bigcup_{\beta<\alpha}V_\beta$, there exist a continuous
function~$f$ defined in~$V$ such that in~$V_\alpha$ the equation
$f(x)=y$ holds.

If also $y\not\in\bigcup_{\beta<\alpha}V_\beta$ then $f$ can be
chosen to be a homeomorphism.
\end{theorem}
\begin{proof}
For the first part of the theorem suppose that we have
$p\in\mathbb{P}_\alpha$ such that $p\Vdash
``\dot{x}\not\in\bigcup_{\beta<\alpha}V_\beta\ \text{and}\
\dot{y}\in\bigcup_{\beta<\alpha}V_\beta"$. We will construct a
fusion sequence below $p$ and define a continuous function~$f$
in~$V$ such that the fusion of the sequence forces that $f(x)=y$
holds in~$V_\alpha$.

Let $p_0=p_1=p$, $n_0=n_1=0$, $F_0=\emptyset$, and choose
$F_1\in[\dom(p)]^{<\omega}$ in such a way that we are building a
fusion sequence. Suppose we have constructed the sequence up to
$i$, we will construct the next element of the fusion sequence.
Let $\{\tau_k:k<K\}$ denote an enumeration of all maps from
$F_{i-1}$ into~$\powerinfront{{n_{i-1}}}$ consistent with~$p_i$.

According to lemma~\ref{lemhlplim} there exists a
$(q_0,j_0)>_{F_i}(p_i\beperk\tau_0,n_i)$ such that for every
$\sigma:F_i\to\powerinfront{{n_i}}$ consistent with~$q_0$ we have
distinct $u_\sigma$'s in $\powerinfront{{m(\tau_0)}}$ (where
$m(\tau_0)$ follows from lemma~\ref{lemhlplim}), such that
$q_0\beperk\sigma\Vdash ``u_\sigma\subset\dot{x}"$. Now use
lemma~\ref{smelt} to construct $r_0\in\mathbb{P}_\alpha$ such that
$r_0\geq q_0$ and $r_0\beperk\tau_0=q_0$.

We now iteratively consider all the $\tau_k$. In the general case
if $\tau_k$ is not consistent with~$r_{k-1}$ then we make sure
that $r_k=r_{k-1}$ and~$j_k=j_{k-1}$. If $\tau_k$ is consistent
with~$r_{k-1}$ we find by lemma~\ref{lemhlplim} a
$(q_k,j_k)>_{F_i}(r_{k-1}\beperk\tau_k,n_i)$ such that for every
$\sigma:F_i\to\powerinfront{{n_i}}$ consistent with~$q_k$ we have
distinct $u_\sigma$'s in $\powerinfront{{m(\tau_k)}}$ such that
$q_k\beperk\sigma\Vdash ``u_\sigma\subset\dot{x}"$. Now use
lemma~\ref{smelt} to construct $r_k\in\mathbb{P}_\alpha$ such that
$r_k\geq r_{k-1}$ and~$r_k\beperk\tau_k=q_k$. After considering
all $\tau_k$ we define $p_{i+1}=r_{K-1}$ and
$n_{i+1}=\max\{j_k:k<K\}$.

If we take a closer look at lemma~\ref{lemhlplim} we can also let
the fusion sequence that we just constructed determine $\dot{y}$.
Because if we have $p\beperk\tau\Vdash ``t_\tau\subset\dot{y}"$,
following the proof of lemma~\ref{lemhlplim} we can make sure that
(by some strengthening of $q^\dagger$ or the $f_\sigma$'s, if
necessary) there exist $t_\sigma$'s in $\powerinfront{{<\omega}}$,
not necessarily distinct, extending~$t_\tau$ such that for
$\sigma:H\to\powerinfront{n}$ consistent with~$q$ we also have
$q\beperk\sigma\Vdash ``t_\sigma\subset\dot{y}"$. So assume we
have done this. We have for every
$\sigma:F_i\to\powerinfront{{n_i}}$ consistent with~$p_{i+1}$
\begin{equation}
\label{equ-in-x-t-in-y} p_{i+1}\beperk\sigma\Vdash
``u_\sigma\subset\dot{x}\ \text{and}\ t_\sigma\subset\dot{y}".
\end{equation}

Now we are ready to define our function~$f$ which will map $x$ in
$V_\alpha$ continuously onto~$y$. Let $f([u_\sigma])\subset
[t_\sigma]$ for all $\sigma:F_i\rightarrow\powerinfront{{n_{i}}}$
and all $i\in\omega$. Then $p_i\beperk\sigma\Vdash ``f(\dot{x})\in
[t_\sigma]"$ for $\sigma:F_i\rightarrow\powerinfront{{n_i}}$
consistent with~$p_i$ and $i\in\omega$. It follows that the
fusion~$q$ forces that in~$V_\alpha$ we have $f(x)=y$. Moreover
$f$ is a continuous function, this follows from
lemma~\ref{lemhlplim}.

For the second part of the theorem suppose that
$p\Vdash``\dot{x},\dot{y}\not\in\bigcup_{\beta<\alpha}V_\beta"$.
Just as in lemma~\ref{lemhlplim} we can choose not only the
$u_\sigma$'s in equation~\ref{equ-in-x-t-in-y} distinct but also
the $t_\sigma$'s for $\sigma\in\Sigma_\tau$ and
$\tau:F_i\to\powerinfront{{n_i}}$ for some $i\in\omega$. With
this, the constructed continuous function~$f$ is actually a
homeomorphism.
\end{proof}
As there are no reals added at limit stages of cofinality larger
than~$\aleph_0$ we have as a corollary to theorems~\ref{thopv} and
\ref{thlim}
\begin{corollary}
For every $\alpha$ and every $\dot{x}$ and $\dot{y}$
$\mathbb{P}_\alpha$-names for reals
in~$V_\alpha\setminus\bigcup_{\beta<\alpha}V_\beta$ there exists a
homeomorphism~$f$ defined in~$V$ such that in $V_\alpha$ we have
$f(x)=y$.
\end{corollary}

\begin{remark}
It is not the case that the $\mathfrak{t}\mathfrak{f}$ number is
the same for all compact metric spaces, e.g. every Cook
continuum~$X$ has $\mathfrak{t}\mathfrak{f}(X)=\mathfrak{c}$ (it
only has the identity and constant mappings as self-maps,
see~\cite{C}). On the other hand, in the Sacks model one has
$\mathfrak{t}\mathfrak{f}(C)=\mathfrak{t}\mathfrak{f}(\mathbb{R})
=\mathfrak{t}\mathfrak{f}([0,1])=\aleph_1$. To see this, observe
that our proof produces, given $x$ and~$y$, two copies of the
Cantor set $A$ and~$B$ containing $x$ and~$y$ respectively and a
continuous map $f:A\to B$, say, with $f(x)=y$. One can then extend
$f$ to a continuous map $\tilde{f}:[0,1]\to [0,1]$ (or
$\tilde{f}:\mathbb{R}\to\mathbb{R}$), whose code will still be
in~$V$.
\end{remark}
\begin{remark}
If $\operatorname{cov}(nowhere\ dense)=\mathfrak{c}$ for the unit
interval $I$, then remark~\ref{remarkcov} shows that
$\mathfrak{t}\mathfrak{f}(I)=\mathfrak{c}$. Suppose that
$\operatorname{cov}(nowhere\ dense)=\kappa<\mathfrak{c}$, for $I$,
then we can cover $I$ by $\kappa$ many Cantor sets
$\{C_\alpha\}_{\alpha<\kappa}$ in such a way that for every two
reals $x$ and~$y$ there exists an $\alpha$ such that $x,y\in
C_\alpha$.  For every $\alpha$ we have a transitive family of
continuous functions $\mathcal{F}_\alpha$ on $C_\alpha$ such that
$\vert\mathcal{F}_\alpha\vert=\mathfrak{t}\mathfrak{f}(C)$. We can
extend every $f\in\mathcal{F}_\alpha$ to a continuous self map
$\tilde{f}$ of~$I$. So $\mathcal{F}=\{\tilde{f}:\ \text{there is
an}\ \alpha<~\kappa\ \text{and}\ f\in\mathcal{F}_\alpha\}$ is a
transitive set of continuous functions on~$I$, and its cardinality
is less than or equal to
$\kappa\times\mathfrak{t}\mathfrak{f}(C)=\mathfrak{t}\mathfrak{f}(C)$.

So if we can cover the unit interval with less than $\mathfrak{c}$
many nowhere dense sets we have
$\mathfrak{t}\mathfrak{f}(I)\leq\mathfrak{t}\mathfrak{f}(C)$.
\end{remark}
%
%
%
%
\section{The cardinal $\mathfrak{t}\mathfrak{f}$ and side-by-side
         Sacks forcing}
In this paper we showed that after adding $\aleph_2$ many Sacks
reals iteratively to a model of $ZFC+CH$ we end up with a model of
$\mathfrak{t}\mathfrak{f}<\mathfrak{c}$. Now consider
$\mathbb{P}\mathbb{S}(\kappa)$, the poset for adding $\kappa$ many
Sacks reals side-by-side (see \cite{B}). We have that
$\mathbb{P}\mathbb{S}(\kappa)$ has the $(2^{\aleph_0})^+$-chain
condition and preserves $\aleph_1$. Suppose that
$\kappa\geq\aleph_1$ and $\cf(\kappa)\geq\aleph_1$. If $V$ is a
model of $CH$ and $G$ is $\mathbb{P}\mathbb{S}(\kappa)$-generic
over $V$, we have in $V[G]$ that $2^{\aleph_0}=\kappa$ and all
cardinals are preserved.

A natural question would be if we get a model of
$\mathfrak{t}\mathfrak{f}<\mathfrak{c}$ when we add $\aleph_2$
many Sacks reals side-by-side to a model of $ZFC+CH$. The answer
to this question is in the negative.

Suppose that $V$ is a model of $ZFC$. Consider the poset
$\mathbb{P}=\mathbb{P}\mathbb{S}(\{1,2,3,4\})$ that adds four
Sacks reals side-by-side to the model $V$. We define
$\mathbb{P}_1$ to be the p.o.-set $\mathbb{P}\mathbb{S}(\{1,2\})$
and $\mathbb{P}_2$ to be the p.o.-set
$\mathbb{P}\mathbb{S}(\{3,4\})$. Suppose $G$ is $\mathbb{P}$
generic over $V$ then $G_{12}=G\beperk\{1,2\}$ is $\mathbb{P}_1$
generic and $G_{34}=G\beperk\{3,4\}$ is $\mathbb{P}_2$ generic
over~$V$. The following holds.
\begin{lemma}
\label{lemsidebyside} In $V[G]$ we have $V[G_{12}]\cap
V[G_{34}]=V$.
\end{lemma}
\begin{proof}
Suppose that $\dot{X}$ is a $\mathbb{P}$ name and $q$ an element
of $\mathbb{P}$ such that $q\Vdash ``\dot{X}\in V[G_{12}]\cap
V[G_{34}]"$. So there exists a $\mathbb{P}_1$ name $\dot{Y}$ and a
$\mathbb{P}_2$ name $\dot{Z}$ such that $q\Vdash
``\dot{X}=\dot{Y}=\dot{Z}"$. Aiming for a contradiction assume
$\dot{X}$ is a name for an object not in~$V$. There exists a
$n\in\omega$ such that $q$ does not decide $n\in\dot{X}$. Now we
have $q_1=q\beperk\{1,2\}$ does not decide $n\in\dot{Y}$, and
$q_2=q\beperk\{3,4\}$ does not decide $n\in\dot{Z}$. So we can
find in $\mathbb{P}_1$ a $r\geq q_1$ such that $r\Vdash
``n\in\dot{Y}"$ and in $\mathbb{P}_2$ a $t\geq q_2$ such that
$t\Vdash ``n\not\in\dot{Z}"$. This gives the contradiction we are
looking for because $r\cup t\Vdash``\dot{Y}\not=\dot{Z}"$ and
$r\cup t\geq q$. So $\dot{X}$ must be a name of an element in~$V$.
\end{proof}
Now we can prove that adding $\aleph_2$ many Sacks reals to a
model of $ZFC+CH$ we do not produce a model of
$\mathfrak{t}\mathfrak{f}<\mathfrak{c}$.
\begin{theorem}
Suppose $V\models CH$ and $G$ is a
$\mathbb{P}\mathbb{S}(\kappa)$-generic filter over~$V$, where
$\kappa\geq\aleph_1$ and $\cf(\kappa)\geq\aleph_1$, then
$V[G]\models\mathfrak{t}\mathfrak{f}=\mathfrak{c}$.
\end{theorem}
\begin{proof}
For every $\alpha<\beta<\kappa$ we have that there exists a
function $f_{\alpha,\beta}\in V[G\beperk\{\alpha,\beta\}]$ mapping
$s_\alpha$ onto~$s_\beta$ or vice versa. This function
$f_{\alpha,\beta}$ is not a member of $V$ for the obvious reason
that assuming that $f_{\alpha,\beta}$ maps $s_\alpha$
onto~$s_\beta$ we get $s_\beta\in V[G\beperk\{\alpha\}]$, which,
of course, is false. Using lemma~\ref{lemsidebyside} and the fact
that $2\kappa=\kappa$ we see that the size of
$\mathfrak{t}\mathfrak{f}$ is at least $\kappa$, because
$f_{2\alpha,2\alpha+1}\not=f_{2\beta,2\beta+1}$ for every
$\alpha\not=\beta$. By theorem~\ref{f<=c<=f+} we are done.
\end{proof}
%
%

\end{document}